\documentclass[11pt]{article}
\renewcommand{\baselinestretch}{1.5}

\usepackage{fancybox,color}
\usepackage{makeidx}
\usepackage{amsmath,amssymb}
\usepackage{verbatim}

\usepackage{exscale}
\usepackage{tabularx}
\usepackage{algorithm,algorithmic}
\usepackage{amsthm}
\usepackage{latexsym}

\usepackage[dvips]{graphicx,epsfig}
\usepackage{url}

\newcommand{\mb}[1]{\mbox{\boldmath $#1$}}

\newenvironment{pf}{\noindent {\bf Proof. }}{\hfill $\square$}

\newcommand{\cX}{{\cal X}}

\newcommand{\cH}{{\cal H}}

\newcommand{\cA}{{\cal A}}

\newcommand{\cQ}{{\cal Q}}
\newcommand{\cB}{{\cal B}}

\newcommand{\bea}{\begin{eqnarray*}}
\newcommand{\eea}{\end{eqnarray*}}

\newcommand{\vect}[1]{{\boldsymbol #1 }}

\newcommand{\inprod}[2]{\langle #1 , #2 \rangle }
\newcommand{\card}[1]{\vert #1 \rvert}
\newcommand{\bc}{\begin{center}}
\newcommand{\ec}{\end{center}}
\newcommand{\bz}{\vect{z}}

\newcommand{\by}{\vect{y}}
\newcommand{\bY}{\vect{Y}}
\newcommand{\refs}[1]{(\ref{#1})}

\newcommand{\bX}{\vect{X}}
\newcommand{\bx}{\vect{x}}

\newcommand{\bZ}{\boldsymbol Z}
\newcommand{\R}{\mathbb R}

\newcommand{\be}{\begin{equation}}
\newcommand{\ee}{\end{equation}}
\newcommand{\beaa}{\begin{eqnarray*}}
\newcommand{\eeaa}{\end{eqnarray*}}
\newcommand{\ben}{\begin{enumerate}}
\newcommand{\een}{\end{enumerate}}
\newcommand{\db}{\hspace*{\fill}{\zapf o}}

\newcommand{\bpn}{\begin{proposition}\twlsf}
\newcommand{\epn}{\db\end{proposition}}
\newcommand{\bdm}{\begin{displaymath}}
\newcommand{\edm}{\end{displaymath}}
\newcommand{\ba}{\begin{array}}
\newcommand{\ea}{\end{array}}

\newcommand{\st}{\mathop{\rm s.t.}}

\newtheorem{lemma}{Lemma}
\newtheorem{proposition}{Proposition}

\newtheorem{remark}{Remark}
\newtheorem{theorem}{Theorem}
\newcommand{\eps}{\epsilon}

\newcommand{\norm}[1]{\| #1 \|}

\def\bV{\vect{V}}
\def\bA{\vect{A}}
\def\lam{\lambda}

\title{A proximal point algorithm for sequential feature extraction applications}
\author{Xuan Vinh Doan
\thanks{Department of Combinatorics and Optimization,
University of Waterloo, 200 University Avenue West, Waterloo, ON N2L 3G1,
Canada, vanxuan@uwaterloo.ca.}
\and Kim-Chuan Toh
\thanks{Department of Mathematics, National University of Singapore,
Blk S17, 10 Lower Kent Ridge Road, Singapore 119076, mattohkc@nus.edu.sg}
\and Stephen Vavasis \thanks{Department of Combinatorics and Optimization,
University of Waterloo, 200 University Avenue West, Waterloo, ON N2L 3G1,
Canada, vavasis@math.uwaterloo.ca.}}

\date{August 2011}

\begin{document}
\maketitle

\begin{abstract}
We propose a proximal point algorithm to solve LAROS problem, that is the problem of finding a ``large approximately rank-one submatrix''. This LAROS problem is used to sequentially extract features in data. We also develop a new stopping criterion for the proximal point algorithm, which is based on the duality conditions of $\eps$-optimal solutions of the LAROS problem, with a theoretical guarantee. We test our algorithm with two image databases and show that we can use the LAROS problem to extract appropriate common features from these images.
\end{abstract}

\section{Introduction}
\label{sec:intro}
Feature extraction is an important application in information retrieval. For example, let us consider a matrix $\mb{A}\in\R_+^{m\times n}$ that represents a database of pixelated and registered grayscale images which have the same size. Each column of $\mb{A}$ corresponds to one image and each row corresponds to a particular pixel position in those images. The value $A_{ij}$ is then the intensity of the $i$th pixel in the $j$th image. A common visual feature represented by the pixels in ${\cal J}\subset\{1,\ldots,n\}$, which occur in a subset of images in ${\cal I}\subset\{1,\ldots,m\}$, can be associated with the approximately rank-one submatrix $\mb{A}({\cal I},{\cal J})$ of the matrix $\mb{A}$. We assume here the features are non-overlapping. If we want to more than one visual feature, we can iteratively find an approximately rank-one submatrix, subtract it from $\mb{A}$ (perhaps modifying the result of the subtraction to ensure that $\mb{A}$ remains nonnegative), and then repeat the procedure. Doan and Vavasis \cite{Doan10} proposed the LAROS problem which tries to find ``large approximately rank-one submatrix''. The proposed convex parametric formulation for the LAROS problem is written as follows:
\be
\label{eq:mprob}
\ba{rl}
\min & \norm{\bX}_\theta \;:=\; \norm{\bX}_*+\theta\norm{\bX}_1\\
\st & \inprod{\mb{A}}{\bX}=1,
\ea
\ee
where $\theta>0$. Here $\norm{\bX}_*$ denotes the nuclear norm of
$\bX$, which is defined to be the sum of the singular values of $\bX$,
and $\norm{\bX}_1$ denotes the sum the absolute values of all the entries of $\bX$. Theoretical properties of LAROS problem have been developed in \cite{Doan10}. In this paper, we investigate algorithms to solve the problem and apply it to find features in data. We will focus on proximal point algorithmic framework, which have recently been studied for nuclear norm minimization (see  Liu et al. \cite{Liu09} and references therein).

Throughout the paper, we use $\norm{\cdot}$ to denote either the
Frobenius norm of a matrix or the Euclidean norm of a vector.
The spectral norm of a matrix $\bX$ is denoted by $\norm{\bX}_2$.

\subsection*{Proximal Point Algorithm}
The proximal point algorithm is based on the \emph{Moreau-Yoshida regularization} of the (non-differentiable) convex optimization problem
\be
\label{eq:prob}
\min_{\bx\in{\cal X}}\phi(\bx),
\ee
where $\cal X$ is a finite-dimensional real Hilbert space and $\phi:{\cal X}\rightarrow (-\infty,\infty]$ is a proper, lower semicontinuous, convex function. For an arbitrary $\lambda>0$, the regularization is defined as
$$
\Phi_{\lambda}(\bx)=\min_{\bz\in{\cal X}}\left(\phi(\bz)+\frac{1}{2\lambda}\norm{\bx-\bz}^2\right),
\quad\forall\,\bx\in{\cal X}.
$$

The above optimization problem has a unique optimal solution $p_{\lambda}(\bx)$ for all $\bx\in {\cal X}$, and $p_{\lambda}$ is called the \emph{proximal point mapping} associated with $\phi$. One of the most important properties of $\Phi_{\lambda}$ and $p_{\lambda}$ is that the set of optimal solutions of \refs{eq:prob} is exactly the set of optimal solutions of the following optimization problem:
\be
\label{eq:rprob}
\min_{\bx\in{\cal X}}\Phi_{\lambda}(\bx),
\ee
where $\Phi_{\lambda}$ is now a \emph{continuously differentiable} convex function defined on $\cal X$ with a globally Lipschitz continuous gradient $\nabla\Phi_{\lambda}$ (with modulus $1/\lambda$). The necessary and sufficient optimality condition of \refs{eq:rprob} can then be expressed as follows:
\be
\label{eq:popt}
\nabla\Phi_{\lambda}(\bx)=\mb{0}\Leftrightarrow p_{\lambda}(\bx)=\bx,
\ee
where $p_{\lambda}$ is a global Lipschitz continuous function with modulus $1$.

The \emph{proximal point} algorithm is an iterative method to solve the problem \refs{eq:prob} that uses the optimality condition written in \refs{eq:popt}. In each iteration, $\bx^{k+1}\approx p_{\lambda_k}(\bx^k)$ according to a sequence $\{\lambda_k\}$ of regularization parameters. The convergence of the algorithm has been studied by Rockafellar \cite{Rockafellar76} in a more general setting of \emph{inclusion} problems with \emph{maximal monotone operators}. Note that the problem \refs{eq:prob} is equivalent to the inclusion problem $\mb{0}\in\partial\phi(\bx)$, where $\partial\phi$ is a maximal monotone operator if $\phi$ is a proper, lower semicontinuous, and convex function. We now ready to study the proximal point mapping for our particular problem. In order to apply the framework, we reformulate Problem \refs{eq:mprob} with a redundant variable as follows:
\be
\label{eq:mmprob}
\ba{rl}
\min & \norm{\bX_1}_*+\theta\norm{\bX_2}_1\\
\st & \inprod{\mb{A}}{\bX_1}=1,\\
\quad & \bX_1 = \bX_2.
\ea
\ee
In addition, to introduce more flexibility into our model, we study Problem (\ref{eq:mmprob}) under the following more general setting:
\begin{eqnarray}
\label{eq:genprob}
\ba{rl}
\min & \norm{\bX_1}_*+\theta\norm{\bX_2}_1\\
\st & \cA(\bX)-b \in \cQ, \quad
\bX = (\bX_1,\bX_2) \in \R^{m\times n}\times\R^{m\times n}
\ea
\end{eqnarray}
where $b\in \cH$,
$\cA: \R^{m\times n}\times\R^{m\times n} \rightarrow \cH$ is a given
linear map, and $\cQ$ is a pointed close convex cone in $\cH$.
Here $\cH$ is a finite-dimensional Hilbert space.
For the problem \refs{eq:mmprob}, we have $\cH = \R\times \R^{m\times n}$,
$\cQ = \{0\}\times \{{\bf 0}\}$,
$b=(1,{\bf 0})$, and $\cA(\bX)=(\inprod{\bA}{\bX_1},\bX_1-\bX_2)$.
Note that the adjoint $\cA^*: \cH \rightarrow \R^{m\times n}\times \R^{m\times n}$
is given by $\cA^* z = (z_1 \bA + \bZ_2, -\bZ_2)$ for any $z=(z_1,\bZ_2) \in \cH$.

%%*********************************************
\section{Primal Proximal Point Algorithm}
We define the function $\phi$ as follows:
\be
\label{eq:phi}
\phi(\bX)=\left\{
\ba{ll}
\norm{\bX_1}_*+\theta\norm{\bX_2}_1, & \bX=(\bX_1,\bX_2)\in{\cal F},\\
+\infty, & \mbox{otherwise},
\ea\right.
\ee
where ${\cal F}$ is the feasible set of the problem \refs{eq:genprob}. The problem \refs{eq:genprob} is then equivalent to the optimization problem
$$
\min_{\bX\in{\cal X}}\phi(\bX),
$$
where ${\cal X}=\R^{m\times n}\times\R^{m\times n}$.

We now introduce dual decision variables $z\in\R^p$ and define the Lagrangian function $L(\bX,z)$,
\be
\label{eq:lagrange}
L(\bX,z)=
\left\{ \begin{array}{ll}
\norm{\bX_1}_*+\theta\norm{\bX_2}_1
+ \inprod{z}{b-\cA(\bX)}
 & \mbox{if $z\in \cQ^*$}
\\[5pt]
-\infty & \mbox{otherwise}
\end{array} \right.
\ee
where $\cQ^*$ is the dual cone of $\cQ$ defined by
$\cQ^* = \{ y \in \cH: \inprod{y}{z} \geq 0, \, \forall \, z\in \cQ\}$. For our problem, $\cQ^*$ is simply the whole space, $\cQ^* = \cH=\R\times\R^{m\times n}$. Clearly, $\displaystyle\phi(\bX)=\sup_{z\in\R^p}L(\bX,z)$. We now calculate the Moreau-Yoshida regularization of $\phi$:
\be
\label{eq:regularized}
\Phi_{\lambda}(\bX)=\min_{\vect{V}\in\cX}
\left(\phi(\mb{V})+\frac{1}{2\lambda}\norm{\bX-\mb{V}}^2\right).
\ee
Applying the strong duality (or minimax theory) result in Rockafellar \cite{Rockafellar70}, we have:
$$
\ba{rcl}
& & \hspace{-0.7cm}
\Phi_{\lambda}(\bX) = \displaystyle\min_{\vect{V}\in\cX} \sup_{z\in\cH}
\left(L(\mb{V},z)+
\frac{1}{2\lambda}\norm{\bX-\mb{V}}^2\right)\\[5pt]
\quad  &=& \displaystyle\sup_{z\in\cH}\min_{\vect{V}\in\cX}
\left(L(\mb{V},z)+\frac{1}{2\lambda}\norm{\bX-\mb{V}}^2\right)
\\[5pt]
\quad  &=& \displaystyle\sup_{z\in\cQ^*}\min_{\vect{V}\in\cX}
\left(\norm{\bV_1}_*+\theta\norm{\bV_2}_1 + \inprod{z}{b-\cA(\bV)}
+\frac{1}{2\lambda}\norm{\bX-\mb{V}}^2\right)
\\[5pt]
&=& \displaystyle\sup_{z\in\cQ^*} \inprod{z}{b}
 + \frac{1}{2\lambda}\norm{\bX}^2 -
 \frac{1}{2\lambda}\norm{\bX + \lambda\cA^*z}^2
 +\min_{\vect{V}\in\cX}
 \left(\norm{\bV_1}_*+\theta\norm{\bV_2}_1
 + \frac{1}{2\lambda}\norm{\bV-(\bX+\lambda\cA^*z)}^2\right)
\ea
$$
Now, consider the first inner minimization problem, we have:
\begin{eqnarray}
&& \hspace{-0.7cm}
\displaystyle\min_{\vect{V}\in\cX}
\left(\norm{\bV_1}_*+\theta\norm{\bV_2}_1
 + \frac{1}{2\lambda}\norm{\bV-(\bX+\lambda\cA^*z)}^2\right)
   \label{eq:1} \\[5pt]
&=&
 \displaystyle\min_{\bV_1}
\left(\norm{\bV_1}_*
 + \frac{1}{2\lambda}\norm{\bV_1-(\bX_1+\lambda \cB_1 z)}^2\right)
 + \theta\displaystyle\min_{\bV_2}
\left(\norm{\bV_2}_1
 + \frac{1}{2\lambda\theta}\norm{\bV_2-(\bX_2+\lambda \cB_2 z)}^2\right)
 \nonumber
\end{eqnarray}
where we have written $\cA^*z = ( \cB_1z, \cB_2z)\in \cX$.
The first optimization problem on the right-hand side is the Moreau-Yoshida
regularization of the nuclear norm function at $\bX_1+\lambda\cB_1z$, and
the problem has an analytical solution given by
\be
\label{eq:nuclearp}
p_{\lambda}^{(1)}(\bX_1+\lambda \cB_1z)
=\mb{U}\mbox{Diag}(\max\{\sigma_i-\lambda,0\})\mb{V}^T,
\ee
which is computable from the singular value decomposition,
$\bX_1+\lambda\cB_1z=\mb{U}\mb{\Sigma}\mb{V}^T$.
In addition, the minimal objective value is given by
$$
\frac{1}{2\lambda}\norm{\bX_1 +\lambda \cB_1z}^2
-\frac{1}{2\lambda}\norm{p_{\lambda}^{(1)}(\bX_1+\lambda \cB_1z)}^2.
$$
Next, we consider the second inner minimization problem on the right-hand side
of (\ref{eq:1}). This optimization problem is the Moreau-Yoshida regularization of
the $l_1$-norm function (with parameter $\lambda\theta$) at
$\bX_2+\lambda\cB_2z$, and it has the following analytical solution:
\begin{eqnarray}
\label{eq:l1p}
\ba{rcl}
p_{\lambda\theta}^{(2)}(\bX_2+\lambda\cB_2z)
=\mbox{sgn}(\bX_2+\lambda\cB_2z)
\circ\max\{\card{\bX_2+\lambda\cB_2z}-\theta\lambda,0\},
\ea
\end{eqnarray}
where $\circ$ is the Hadamard product (or entrywise product) and
$\mbox{sgn}$ is the (entrywise) sign function. The corresponding
minimal objective value is given by
$$
\frac{1}{2\lambda\theta}\norm{\bX_2 +\lambda\cB_2z}^2
-\frac{1}{2\lambda\theta}
\norm{p_{\lambda\theta}^{(2)}(\bX_2+\lambda\cB_2z)}^2.
$$

Combining these two results, we can compute $\Phi_{\lambda}(\bX)$ as follows:
\be
\label{eq:fregularized}
\Phi_{\lambda}(\bX)=\frac{1}{2\lambda}\norm{\bX}^2
+\sup_{z\in \cQ^*}\left(\inprod{z}{b}
-\frac{1}{2\lambda}
   \norm{p_{\lambda}^{(1)}(\bX_1+\lambda \cB_1z)}^2
-\frac{1}{2\lambda}
   \norm{p_{\lambda\theta}^{(2)}(\bX_2+\lambda\cB_2z)}^2\right),
\ee
where $p_{\lambda}^{(1)}$ and $p_{\lambda\theta}^{(2)}$ are defined in
\refs{eq:nuclearp} and \refs{eq:l1p} respectively. Now define
\be
\label{eq:innerprob}
\Theta_{\lambda}(\bX,z)
= \inprod{z}{b}
-\frac{1}{2\lambda}
   \norm{p_{\lambda}^{(1)}(\bX_1+\lambda\cB_1z)}^2
-\frac{1}{2\lambda}
   \norm{p_{\lambda\theta}^{(2)}(\bX_2+\lambda\cB_2z)}^2
\ee
and consider
$$
z_{\lambda}(\bX)\in\arg\sup_{z\in\cQ^*}\Theta_{\lambda}(\bX,z).
$$
Applying the saddle point theorem in Rockafellar \cite{Rockafellar70},
we obtain the proximal point mapping associated with $\phi$ as follows:
\be
\label{eq:prox}
p_{\lambda}(\bX)=
\left(p_{\lambda}^{(1)}(\bX_1+\lambda\cB_1z_\lambda(\bX)),\,
p_{\lambda\theta}^{(2)}(\bX_2+\lambda\cB_2 z_\lambda(\bX)) \right).
\ee
The primal proximal point algorithm (primal PPA) has the following
template.

\vspace{10pt}

\fbox{\parbox{0.95\textwidth}{\noindent{\bf The Primal PPA.} Given $\bX^0\in \cX$, $\lambda_0 > 0$ and $\varepsilon > 0$, perform
the following loop:
\begin{description}
\item[Step 1.] Find an (approximate) optimal solution
\begin{eqnarray}
 z^{k}\in\arg\sup_{z\in\cQ^*} \Theta_{\lambda_k}(\bX^k,z),
\end{eqnarray}
     where $\Theta_{\lambda_k}$ is defined in \refs{eq:innerprob}.
\item[Step 2.] Update
\begin{eqnarray}
\bX^{k+1}_1= p_{\lambda_k}^{(1)}(\bX^k_1+\lambda \cB_1z^k),
\quad
\bX^{k+1}_2= p_{\lambda_k\theta}^{(2)}(\bX^k_2+\lambda\cB_2z^k)
\end{eqnarray}
according to the proximal point mapping in \refs{eq:prox}.
\item[Step 3.] If $\norm{\bX^{k+1}-\bX^k}/\lambda_k < \varepsilon$, stop;
else, update $\lambda_k$, end
\end{description}}}

\section{Dual Proximal Point Algorithm}

The dual problem associated with \refs{eq:genprob} is given as follows:
\begin{eqnarray}
\max_{y\in\cH} g(y)
\end{eqnarray}
where $g$ is the concave function defined by
\begin{eqnarray}
  g(y) = \left\{ \begin{array}{ll}
  \inf \{ \norm{\bX_1}_* +\theta\norm{\bX_2}_1
   + \inprod{y}{b-\cA(\bX)} \,: \,
  \bX = (\bX_1,\bX_2) \in \cX
  \} & \mbox{if $y\in\cQ^*$}
   \\[5pt]
   -\infty &\mbox{otherwise}
   \end{array} \right.
\end{eqnarray}
The Moreau-Yoshida regularization of $g$ is given by
\begin{eqnarray}
 G_\lambda(y) &:=&
  \max_{z\in \cH} \{g(z) - \frac{1}{2\lambda}\norm{z-y}^2\}
  \nonumber \\[5pt]
 &=&\max_{z\in \cQ^*}
 \inf_{\bX\in\cX} \{ \norm{\bX_1}_* +\theta\norm{\bX_2}_1
   + \inprod{z}{b-\cA(\bX)}
   - \frac{1}{2\lambda}\norm{z-y}^2\}
   \nonumber \\[5pt]
 &=&
 \inf_{\bX\in\cX}
 \max_{z\in \cQ^*} \{ \norm{\bX_1}_* +\theta\norm{\bX_2}_1
   + \inprod{z}{b-\cA(\bX)}
   - \frac{1}{2\lambda}\norm{z-y}^2\}
   \nonumber \\[5pt]
 &=&  -\frac{1}{2\lambda}\norm{y}^2 +
 \inf_{\bX\in\cX}
  \Big\{ \norm{\bX_1}_* +\theta\norm{\bX_2}_1
   + \Psi_\lam(\bX;y)\Big\}
\end{eqnarray}
where
\begin{eqnarray}
\label{eq:Psi}
 \Psi_\lam(\bX;y) =
 \frac{1}{2\lambda} \norm{\Pi_{\cQ^*}(y+\lambda(b-\cA(\bX)))}^2.
\end{eqnarray}
Note that
$\nabla_\bX \Psi_\lam(\bX;y)=
-\cA^* \Pi_{\cQ^*}(y+\lambda(b-\cA(\bX)))$. The dual algorithm can then be written as follows.

\vspace{10pt}

\fbox{\parbox{0.95\textwidth}{\noindent{\bf The Dual PPA.} Given a
tolerance $\varepsilon>0$. Input $y^0\in \cQ^*$ and $\lam_0>0$. Set
$k:=0$. Iterate:
\begin{description}
\item [Step 1.] Find an approximate minimizer
   \begin{equation}
     \label{eq-inner-dual}
   \bX^{k}\approx\arg\inf_{X\in\cX}
  \Big\{\norm{\bX_1}_* + \theta\norm{\bX_2}_1
  +\Psi_{\lam_k}(\bX;y^k)\Big\},
  \end{equation}
  where $\Psi_{\lam_k}(\bX; y^k)$ is defined as in \refs{eq:Psi}.
\item [Step 2.] Compute
  \be
  y^{k+1} \;=\; \Pi_{{\cal Q}^*}\big[y^k+\lam_k(b-\cA(\bX^{k}))\big].
  \ee
\item [Step 3.] If
$\norm{(y^{k}-y^{k+1})/\lam_k}\leq\varepsilon$; stop; else; update
$\lam_k$ ; end.
\end{description}}}

\section{Implementation Issues}

\subsection{Primal proximal point algorithm}

For the primal PPA, the most important issue we have to
address first is how to solve the inner problem
$\displaystyle\sup_{z\in \cQ^*} \Theta_{\lambda}(\bX,z)$.
We have that $\Theta_{\lambda}$ is a concave function in $z$
due to the linearity in $z$ of the Lagrangian function $L$.
From the general gradient formulation
$\displaystyle\nabla\Phi_{\lambda}(\bx)=
\frac{1}{\lambda}\left(\bx-p_{\lambda}(\bx)\right)$, we have
that
$
\nabla\norm{p_{\lambda}^{(i)}(\bX_i)}^2=
p_{\lambda}^{(i)}(\bX_i), \; i=1,2.
$
Thus $\Theta_{\lambda}$ is continuously differentiable with
\begin{eqnarray}
 \nabla_z \Theta_{\lambda}(\bX,z)
  &=& b - \cB_1^* p_{\lambda}^{(1)}(\bX_1+\lambda \cB_1z)
   - \cB_2^* p_{\lambda\theta}^{(2)}(\bX_2+\lambda \cB_2z).
\end{eqnarray}
Note that for the problem \refs{eq:mmprob}, we have
$\cB_1 z = \bA z_1 + \bZ_2$, $\cB_2 z = -\bZ_2$ for
$z=(z_1,\bZ_2) \in \R\times \R^{m\times n}$.
Correspondingly, we have
$\cB_1^* (\bX_1) = (\inprod{\bA}{\bX_1},\bX_1)$
and $\cB_2^* (\bX_2) = (0,-\bX_2)$ for any $\bX_1,\bX_2\in \R^{m\times n}$
and
\begin{eqnarray}
 \nabla_z \Theta_{\lambda}(\bX,z)
  =\Big(1-\inprod{\bA}{ p_{\lambda}^{(1)}(\bX_1+\lambda (\bA z_1 + \bZ_2))},\,
   p_{\lambda\theta}^{(2)}(\bX_2-\lambda\bZ_2 )-
   p_{\lambda}^{(1)}(\bX_1+\lambda (\bA z_1 + \bZ_2))\Big).
\end{eqnarray}
In addition, using the global Lipschitz continuity (with modulus $1$) of two
proximal point mappings, $p_{\lambda}^{(1)}$ and $p_{\lambda\theta}^{(2)}$,
we can show that the gradient $\nabla_z\Theta_{\lambda}$ is globally
Lipschitz continuous with modulus $\lambda(\norm{\mb{A}}_2^2+2)$.

With all these properties of $\Theta_{\gamma}$, we can solve the inner problem
using \emph{first-order gradient-based} methods such as steepest descent method.

The second issue is that these inner problems are typically only solved approximately
which results in inexact proximal point mappings. For inexact proximal point method, Rockafellar \cite{Rockafellar76} provides two convergence criteria for global and local convergence. Based on the aforementioned convergence criteria,
Liu et al. \cite{Liu09} have proposed some checkable stopping criteria for the inner problems to maintain global (and local) convergence of the proposed (inexact) proximal point method (for nuclear norm minimization problems). We can extend these stopping criteria for our problem.

The third issues is to calculate a partial singular value decomposition in order to compute the proximal point mapping of the nuclear norm function (the computation of the proximal point mapping of the $l_1$-norm function is straightforward). As in Liu et al. \cite{Liu09}, we use a Lanczos bidiagonalization algorithm with partial reorthogonalization to compute a partial singular value decomposition. We also need heuristics to set the number of singular values required to be computed with this algorithm.

\subsection{Dual proximal point algorithm}

For the dual PPA, we need to look for a method to solve the inner problem $\displaystyle \inf_{X\in\cX}
  \Big\{\norm{\bX_1}_* + \theta\norm{\bX_2}_1
  +\Psi_{\lam_k}(\bX;y^k)\Big\}$. Similar to the approach proposed in Liu et al. \cite{Liu09}, we will apply the accelerated proximal gradient algorithm
  \cite{BT09} for this problem. According to Toh and Yun \cite{Toh09}, we solve the problem $\displaystyle \min_{\bX} P(\bX)+f(\bX)$, where $P(\bX)=\norm{\bX_1}_*+\theta\norm{\bX_2}_1$ and $f(\bX)=\Psi_\lam(\bX;y)$. We have  that the gradient $\nabla_\bX \Psi_\lam(\bX;y)=
-\cA^* \Pi_{\cQ^*}(y+\lambda(b-\cA(\bX)))$ is globally Lipschitz continuous with modulus $L=\lambda(\norm{\mb{A}}_2^2+2)$. The proximal gradient algorithm in each iteration needs to solve the following quadratic approximation of the sum $P(\bX)+f(\bX)$ at the current solution $\bY$:
\begin{align*}
Q_t(\bX;\bY)&=P(\bX)+f(\bY)+\inprod{\nabla f(\bY}{\bX-\bY}+\frac{t}{2}\norm{\bX-\bY}_F^2\\
\quad & =P(\bX)+\frac{t}{2}\norm{\bX-G_t(\bY)}_F^2 + f(\bY)-\frac{1}{2t}\norm{\nabla f(\bY)}_F^2,
\end{align*}
where $\displaystyle G_t(\bY)=\bY-\frac{1}{t}\nabla f(\bY)$. This function is a strongly convex function in $\bX$ and hence it has a unique minimizer $S_t(\bY)$. We have that
\begin{align*}
P(\bX)+\frac{t}{2}\norm{\bX-G_t(\bY)}_F^2 &=\norm{\bX_1}_*+\theta\norm{\bX_2}_1+\frac{t}{2}\left(\norm{\bX_1-G_t^1(\bY)}_F^2+\norm{\bX_2-G_t^2(\bY)}_F^2\right)\\
\quad & = \left(\norm{\bX_1}_*+\frac{t}{2}\norm{\bX_1-G_t^1(\bY)}_F^2\right)+\left(\theta\norm{\bX_2}_1+\frac{t}{2}\norm{\bX_2-G_t^2(\bY)}_F^2\right),
\end{align*}
where $G_t(\bY)=(G_t^1(\bY),G_t^2(\bY))$. Thus the minimizer $S_t(\bY)=(S_t^1(\bY),S_t^2(\bY))$, where $S_t^1(\bY)$ is the minimizer of the problem $\displaystyle\min_{\bX_1}\left(\norm{\bX_1}_*+\frac{t}{2}\norm{\bX_1-G_t^1(\bY)}_F^2\right)$ and $S_t^2(\bY)$ is the minimizer of the optimization problem $\displaystyle\min_{\bX_2}\theta\norm{\bX_2}_1+\frac{t}{2}\norm{\bX_2-G_t^2(\bY)}_F^2$. Similar to the previous section, the analytical solutions for these two optimization problems can be calculated and they are:
\begin{eqnarray}
S_t^1(\bY)=p_{t^{-1}}^{(1)}(G_t^1(\bY)),\quad S_t^2(\bY)=p_{t^{-1}\theta}^{(2)}(G_t^2(\bY)).
\label{eq-S}
\end{eqnarray}
Finally, the proximal gradient algorithm for our problem can be described as follows. Given $\tau_0=\tau_{-1}=1$ and $\bX^0=\bX^{-1}$, each iteration includes the following steps
\ben
\item Calculate $\displaystyle\bY^k=\bX^k+\frac{\tau_{k-1}-1}{\tau_k}\left(\bX^{k}-\bX^{k-1}\right)$
\item Update $\bX^{k+1}=S_{t^k}(\bY^k)$ according the formulas in
   (\ref{eq-S}).
\item Update $\displaystyle \tau_{k+1}=\frac{1}{2}\left(\sqrt{1+4\tau_k^2}+1\right)$
\een

The update of $\tau_k$ in the third step is to make sure that $\tau_{k+1}^2-\tau_{k+1}\leq \tau_k^2$ and $\tau_{k+1}\geq 1$, a convergence condition of the proximal gradient algorithm. We also need to have the update rule for the remaining parameter $t_k$, which affects the quadratic approximation of the function $f$ at $\bY$. Since the gradient $\nabla f(\bY)$ is Lipschitz continuous with modulus $L$, for all $t\geq L$, we have:
$$
P(S_t(\bY))+f(S_t(\bY))\leq Q_t(S_t(\bY);\bY).
$$
In order to have a better approximation, we would like to have smaller $t$ and in the accelerated proximal gradient algorithm, we will use line search to find $t_k<L$ such that the above condition is still satisfied, starting with $t_1=L$. More details can be found in Toh and Yun \cite{Toh09}.

\section{Sparse Structure of Rank-One Optimal Solutions}
The proposed proximal point algorithm is a first-order iterative method, which normally does not have fast convergence. Applying duality results obtained by Doan and Vavasis \cite{Doan10} for Problem \refs{eq:mprob}, we would like to study better stopping criteria for the proposed proximal point algorithms. We focus on the case when Problem \refs{eq:mprob} has a rank-one optimal solution $\bX=\sigma\mb{u}\mb{v}^T$ since rank-one optimal solutions are what we are looking for in general. The purpose of the termination test is to obtain the correct supports of $\mb{u}$ and $\mb{v}$, that is, the positions of their nonzero entries with a guarantee certificate when we only have approximate values for $\mb{u}$ and $\mb{v}$ from the proposed first-order algorithm. Although the technique in this section is developed for Problem \refs{eq:mprob}, similar ideas can be applied to other proposed formulations with the nuclear norm such as the matrix completion problem. In particular, a test like this for the matrix completion problem can be used to rigorously establish the correct rank of the optimal solution from approximate solutions obtained from a first-order method.

Now let us consider the rank-one optimal solution $\bX$ of the following form
$$
\bX=\begin{pmatrix}
\sigma_1\mb{u}_1\mb{v}_1^T & \mb{0}\\
\mb{0} & \mb{0}
\end{pmatrix},
$$
where $\mb{u}_1\geq\mb{0}$ is a unit vector in $\R^M$, $M\leq m$, and $\mb{v}_1\geq\mb{0}$ is a unit vector in $\R^N$, $N\leq n$. If $\mb{u}_1$ and $\mb{v}_1$ are determined, $\sigma_1$ can be easily calculated to satisfy the optimality condition $\norm{\bX}_{\theta}=1$ (we assume here $\mb{A}\neq\mb{0}$). Note that in general, the rank-one optimal solution $\bX$ could have a different block structure. However, without loss of generality, we can assume that $\mb{u}_1\mb{v}_1^T$ forms an upper left principal submatrix of $\bX$ for ease of exposition. Under this assumption, we can set $\mb{u}=[\mb{u}_1;\mb{0}]\in\R^m$ and $\mb{v}=[\mb{v}_1;\mb{0}]\in\R^n$ with $\sigma=\sigma_1$. Similar to Theorem 5 in Doan and Vavasis \cite{Doan10}, we can then write the optimality conditions as follows:
\begin{quote}
There exists $\mb{W}\in \R^{m\times n}$ and $\mb{V} \in
\R^{m\times n}$ such that
\begin{eqnarray}
& \mb{A} \;=\; \norm{\mb{A}}_{\theta}^*(\mb{u}\mb{v}^T+\mb{W}) +  \theta\norm{\mb{A}}_{\theta}^*\mb{V} &
\\[5pt]
\label{eq-KKT}
&
\norm{\mb{W}}_2\leq 1, \quad \mb{W}^T\mb{u}=\mb{0}, \quad \mb{W}\mb{v}=\mb{0}, \quad
\mb{V}_{11}=\mb{E}_{M\times N}, \quad
\norm{\mb{V}}_{\infty}\leq 1,
& \nonumber
\end{eqnarray}
where $\mb{E}_{M \times N}$ is the $M\times N$ matrix of all ones.
\end{quote}

Letting $\lambda = 1/\norm{\mb{A}}_{\theta}^*$ and splitting all matrices into four subblocks according to the sparse structure of $\bX$, we obtain the following detailed optimality conditions:
\begin{eqnarray}
&\mbox{$(1,\mb{u}_1,\mb{v}_1)$ is a singular triple of $\lambda\mb{A}_{11}-\theta\mb{V}_{11}$, and $\mb{W}_{11}=(\lambda\mb{A}_{11}-\theta\mb{V}_{11})-\mb{u}_1
\mb{v}_1^T$}& \label{eq-KKT-1}
\\
&\mbox{$\mb{W}_{12}=\lambda\mb{A}_{12}- \theta\mb{V}_{12}$, $\mb{W}_{12}^T\mb{u}_1=\mb{0}$, and $\norm{\mb{V}_{12}}_{\infty}\leq 1$} & \label{eq-KKT-2}
\\
&\mbox{$\mb{W}_{21}=\lambda\mb{A}_{21}- \theta\mb{V}_{21}$, $\mb{W}_{21}\mb{v}_1=\mb{0}$, and $\norm{\mb{V}_{21}}_{\infty}\leq 1$}& \label{eq-KKT-3}
\\
&\mbox{$\mb{W}_{22}=\lambda\mb{A}_{22}- \theta\mb{V}_{22}$, and $\norm{\mb{V}_{22}}_{\infty}\leq 1$}& \label{eq-KKT-4}
\\
&\mbox{$\norm{\mb{W}}_2 \leq 1$}.&
\label{eq-KKT-5}
\end{eqnarray}

The following lemma shows how to find $\mb{u}_1$, $\mb{v}_1$ and $\lambda$ (or $\norm{\mb{A}}_{\theta}^*$) from the first optimality condition.
\begin{lemma}
\label{lem:newton}
If $(\lambda,\mb{u}_1,\mb{v}_1)$ satisfies \refs{eq-KKT-1}, then $\bx = (\lambda,\mb{u}_1,\mb{v}_1)$ is a solution of the following system of nonlinear equations
\be
\label{eq:newton}
P(\bx)=\begin{pmatrix}
(\lambda\mb{A}_{11}-\theta\mb{V}_{11})\mb{v}_1-\mb{u}_1\\
(\lambda\mb{A}_{11}-\theta\mb{V}_{11})^T\mb{u}_1-\mb{v}_1\\
\mb{u}_1^T\mb{u}_1 - 1
\end{pmatrix} = \mb{0}.
\ee
\end{lemma}

\begin{pf}
It is easily to see that $\mb{v}_1^T\mb{v}_1=\mb{v}_1^T(\lambda\mb{A}_{11}-\theta\mb{V}_{11})^T\mb{u}_1=\mb{u}_1^T\mb{u}_1=1$ and the first two equations indicate that $(1,\mb{u}_1,\mb{v}_1)$ is a singular triple of $\lambda\mb{A}_{11}-\theta\mb{V}_{11}$.
\end{pf}

The system of equations in \refs{eq:newton} has $M+N+1$ variables and $M+N+1$ equations, which can be solved using Newton method. One of the convergence results of the Newton's method is the Kantorovich theorem, which is given as follows (see Tapia \cite{Tapia71}).

\begin{theorem}[Kantorovich]
\label{thm:kantorovich}
Assume that $P$ is defined and is Fr\'echet differentiable at each point in a given open convex set $D_0$ and for some $\bx_0\in D_0$ that $[P'(\bx_0)]^{-1}$ exists and that
\ben
\item[(i)] $\norm{[P'(\bx_0)]^{-1}}\leq B$,
\item[(ii)] $\norm{[P'(\bx_0)]^{-1}P(\bx_0)}\leq \eta$, and
\item[(iii)] $\norm{P'(\bx)-P'(\by)}\leq K\norm{\bx-\by}$, for all $\bx$ and $\by$ in $D_0$,
\een
with $\displaystyle h=BK\eta\leq\frac{1}{2}$.

Let $\Omega_*=\left\{\bx\,|\,\norm{\bx-\bx_0}\leq t^*\right\}$, where $\displaystyle t^*=\left(\frac{1-\sqrt{1-2h}}{h}\right)\eta$. Now if $\Omega_*\subset D_0$, then the Newton iterates, $\bx_{k+1}=\bx_k-[P'(\bx_k)]^{-1}P(\bx_k)$, are well defined, remain in $\Omega_*$, and converge to $\bx^*\in\Omega_*$ such that $P(\bx^*)=\mb{0}$. In addition,
$$
\norm{\bx^*-\bx_k}\leq\frac{\eta}{h}\left[\frac{\left(1-\sqrt{1-2h}\right)^{2^k}}{2^k}\right],\quad k=0,1,2,\ldots.
$$
\end{theorem}

According to Theorem \ref{thm:kantorovich}, if we can find $\bx_0$ with the corresponding parameter $h\leq 1/2$, then for an arbitrary $\eps>0$, an $\eps$-solution $\bx$ such that $\norm{\bx-\bx^*}\leq\eps$, can be achieved after a finite number of Newton iterations. Now assuming that we have obtained an $\eps$-solution $(\lambda,\mb{u}_1,\mb{v}_1)$ of the system of equations in \refs{eq:newton}, we would like to characterize the sufficient conditions which guarantee that the corresponding solution $(\lambda^*,\mb{u}_1^*,\mb{v}_1^*)$ defines the optimal solution $\bX$ of \refs{eq:mprob} as described above. The following proposition shows these sufficient conditions.

\begin{proposition}
\label{prop:eps}
Consider an $\eps$-solution $(\lambda,\mb{u}_1,\mb{v}_1)$ of the system of equations in \refs{eq:newton}, $0<\eps<1/2$. The corresponding solution $(\lambda^*,\mb{u}_1^*,\mb{v}_1^*$) defines the rank-one optimal solution $\bX^*$,
$$
\bX^*=\begin{pmatrix}
\sigma_1^*\mb{u}_1^*(\mb{v}_1^*)^T & \mb{0}\\
\mb{0} & \mb{0}
\end{pmatrix},
$$
of \refs{eq:mprob} if there exist $\mb{W}$ and $\mb{V}$ that satisfy the following conditions:
\ben
\item[(i)] $\mb{W}_{11}=(\lambda\mb{A}_{11}-\theta\mb{V}_{11})-\mb{u}_1\mb{v}_1^T$ and $\mb{V}_{11}=\mb{E}_{M\times N}$,
\item[(ii)] $\mb{W}_{12}=\lambda\mb{A}_{12}-\theta\mb{V}_{12}$, $\mb{W}_{12}^T\mb{u}_1=\mb{0}$, and $\norm{\mb{V}_{12}}_{\infty}\leq 1-\theta^{-1}(\norm{\mb{A}_{12}}_{\infty}+5)\eps$,
\item[(iii)] $\mb{W}_{21}=\lambda\mb{A}_{21}-\theta\mb{V}_{21}$, $\mb{W}_{21}\mb{v}_1=\mb{0}$, and $\norm{\mb{V}_{21}}_{\infty}\leq 1-\theta^{-1}(\norm{\mb{A}_{21}}_{\infty}+5)\eps$,
\item[(iv)] $\mb{W}_{22}=\lambda\mb{A}_{22}-\theta\mb{V}_{22}$, and $\norm{\mb{V}_{22}}_{\infty}\leq 1$, and
\item[(v)] $\norm{\mb{W}}_2\leq 1-(\norm{\mb{A}}_2+7.5)\eps$.
\een
\end{proposition}

\begin{remark}
In order to test stopping conditions specified in Proposition \ref{prop:eps}, we need to start with an $\eps$-approximate of the optimal solution  $(\lambda^*,\mb{u}_1^*,\mb{v}_1^*)$, where $\mb{u}_1^*$ and $\mb{v}_1^*$ are unit vectors. It is therefore better to solve the problem where $\lambda^*=1/\norm{\mb{A}}_{\theta}^*$ has the same magnitude as entries of $\mb{u}_1^*$ and $\mb{v}_1^*$. Heuristically, we could scale $\mb{A}$ so that $\norm{\mb{A}}_2=1$ to (partially) control the magnitude of $\lambda^*$.
\end{remark}

\begin{pf} Suppose we are given $\mb{W}$ and $\mb{V}$ which satisfy
the conditions (i)--(v).
We will construct $\mb{W}^*$ and $\mb{V}^*$ from $\mb{W}$ and $\mb{V}$ and prove that they satisfy all optimality conditions in
(\ref{eq-KKT-1})--(\ref{eq-KKT-5}) when combining with the solution $(\lambda^*,\mb{u}_1^*,\mb{v}_1^*)$ of (\ref{eq:newton}). We start with the $(1,1)$ subblock . Clearly, we need $\mb{V}_{11}^*=\mb{V}_{11}=\mb{E}_{M\times N}$ and $\mb{W}_{11}^* = (\lambda^*\mb{A}_{11}-\theta\mb{V}_{11}^*)-\mb{u}_1^*(\mb{v}_1^*)^T$. We have:
$$
\mb{W}_{11}^*-\mb{W}_{11}=(\lambda^*-\lambda)\mb{A}_{11}-(\mb{u}_1^*(\mb{v}_1^*)^T-\mb{u}_1\mb{v}_1^T).
$$
Since $(\lambda,\mb{u}_1,\mb{v}_1)$ is an $\eps$-solution, we have that $\max\left\{\card{\Delta\lambda},\norm{\Delta\mb{u}_1},
\norm{\Delta\mb{v}_1}\right\}\leq\eps$, where $\Delta\lambda=\lambda-\lambda^*$, $\Delta\mb{u}_1=\mb{u}_1-\mb{u}_1^*$, and $\Delta\mb{v}_1=\mb{v}_1-\mb{v}_1^*$. Hence
$$
\ba{rl}
\norm{\mb{u}_1^*(\mb{v}_1^*)^T-\mb{u}_1\mb{v}_1^T}
&=\norm{\mb{u}_1^*(\mb{v}_1^*)^T-(\mb{u}_1^*+\Delta\mb{u}_1)(\mb{v}_1^*+\Delta\mb{v}_1)^T}\\
\quad & = \norm{\Delta\mb{u}_1(\mb{v}_1^*)^T+\mb{u}_1^*\Delta\mb{v}_1^T
+\Delta\mb{u}_1\Delta\mb{v}_1^T} \; \leq \; 2\eps + \eps^2,
\ea
$$
since $\norm{\mb{u}_1^*}=\norm{\mb{v}_1^*}=1$.
%Thus we have:
%$$
%\norm{\mb{W}_{11}^*-\mb{W}_{11}}\leq \eps(\norm{\mb{A}_{11}}+2)+\eps^2.
%$$

We continue with the $(2,2)$ subblock. Let $\mb{V}_{22}^*=\mb{V}_{22}$ and $\mb{W}_{22}^*=\lambda^*\mb{A}_{22}-\theta\mb{V}_{22}^*$, we have:
$$
\mb{W}_{22}^*-\mb{W}_{22} = (\lambda^*-\lambda)\mb{A}_{22}.
$$
Now consider the $(1,2)$ subblock, we would like to construct $\mb{W}_{12}^*$  that is close to $\mb{W}_{12}$ and satisfies the condition that $(\mb{W}_{12}^*)^T\mb{u}_1^*=\mb{0}$. We will use appropriate Householder matrices to construct $\mb{W}_{12}^*$ as follows.
For two different unit vectors $\bx$ and $\by$, the Householder matrix $\mb{Q}=\mb{I}-2\bz\bz^T$ with $\displaystyle\bz=\pm\frac{\bx-\by}{\norm{\bx-\by}}$ transforms $\bx$ to $\by$ and vice versa. In other words, $\mb{Q}\bx=\by$ and $\mb{Q}\by=\bx$. The Householder matrix $\mb{Q}$ is symmetric and orthonormal. Now  consider $\bar{\mb{u}}_1=\mb{u}_1/\norm{\mb{u}_1}$. Note that since $\norm{\mb{u}_1^*}=1$ and $\norm{\Delta\mb{u}_1}\leq\eps$, we have
that
 $\card{\norm{\mb{u}_1}-1}\leq\eps$, which implies $\norm{\Delta\bar{\mb{u}}_1}\leq 2\eps$, where $\Delta\bar{\mb{u}}_1=\bar{\mb{u}}_1-\mb{u}_1^*$. We define $\bx=\displaystyle-\frac{\mb{u}_1^*+\bar{\mb{u}}_1}{\norm{\mb{u}_1^*+\bar{\mb{u}}_1}}$ and consider two Householder matrices, $\mb{Q}_1$ and $\mb{Q}_2$, which transform $\bar{\mb{u}}_1$ to $\bx$ and $\bx$ to $\mb{u}_1^*$, respectively. Let us define
$$
\mb{w}_1=\bar{\mb{u}}_1+\frac{\mb{u}_1^*+\bar{\mb{u}}_1}{\norm{\mb{u}_1^*+\bar{\mb{u}}_1}},\quad\mb{w}_2 = \mb{u}_1^* + \frac{\mb{u}_1^*+\bar{\mb{u}}_1}{\norm{\mb{u}_1^*
+\bar{\mb{u}}_1}},
$$
then $\mb{Q}_1$ and $\mb{Q}_2$ can be constructed with $\bz_1=\mb{w}_1/\norm{\mb{w}_1}$ and $\bz_2=\mb{w}_2/\norm{\mb{w}_2}$, respectively. We have
$$
\ba{rl}
\mb{w}_1^T\mb{w}_1 & =\displaystyle\left(\bar{\mb{u}}_1+\frac{\mb{u}_1^*
+\bar{\mb{u}}_1}{\norm{\mb{u}_1^*+\bar{\mb{u}}_1}}\right)^T
\left(\bar{\mb{u}}_1+\frac{\mb{u}_1^*+\bar{\mb{u}}_1}
{\norm{\mb{u}_1^*+\bar{\mb{u}}_1}}\right)
\; =\; 2 + \norm{\mb{u}_1^*+\bar{\mb{u}}_1},
\ea
$$
or $\norm{\mb{w}_1}=\sqrt{2 + \norm{\mb{u}_1^*+\bar{\mb{u}}_1}}$. Similarly, we can also show that $\norm{\mb{w}_2}=\norm{\mb{w}_1}=\sqrt{2 + \norm{\mb{u}_1^*+\bar{\mb{u}}_1}}$. Thus we have
$$
\Delta\bz_1=\bz_1-\bz_2 =\frac{1}{\sqrt{2 + \norm{\mb{u}_1^*+\bar{\mb{u}}_1}}}\Delta\bar{\mb{u}}_1.
$$
Hence $\displaystyle\norm{\Delta\bz_1} \leq \frac{1}{\sqrt{4 -\norm{\Delta\bar{\mb{u}}_1}}}\norm{\Delta\bar{\mb{u}}_1}
\leq \frac{2}{\sqrt{3}}\eps$.

Now consider $\mb{Q}_{12}=\mb{Q}_{2}\mb{Q}_{1}$, we have: $\mb{Q}\bar{\mb{u}}_1=\mb{u}_1^*$ and $\mb{Q}$ is also an orthonormal matrix. Define $\mb{W}_{12}^*=\mb{Q}_{12}\mb{W}_{12}$, we have: $\mb{W}_{12}^*$ satisfies the condition $(\mb{W}_{12}^*)^T\mb{u}_1^*=\mb{0}$ since $\mb{W}_{12}^T\bar{\mb{u}}_1=\mb{0}$. We can then select $\mb{V}_{12}^*=(\lambda^*\mb{A}_{12}-\mb{W}_{12}^*)/\theta$.
Thus
$$
\mb{W}_{12}^*-\mb{W}_{12}=(\mb{Q}_{12}-\mb{I})\mb{W}_{12},\quad\mb{V}_{12}^*-\mb{V}_{12}=\frac{1}{\theta}\left[(\lambda^*-\lambda)\mb{A}_{12}-(\mb{Q}_{12}-\mb{I})\mb{W}_{12}\right].
$$
We have
$$
\mb{Q}_{12}-\mb{I} =\left(\mb{I}-2\mb{z}_2\mb{z}_2^T\right)\left(\mb{I}-2\mb{z}_1\mb{z}_1^T\right)-\mb{I}\\
\; =\; 2\bz_2\Delta\bz_1^T- 2\Delta\bz_1 \bz_1^T + 4 (\bz_2^T\Delta\bz_1)\bz_2\bz_1^T.
$$
Thus
$$
\frac{1}{4}\norm{\mb{Q}_{12}-\mb{I}}^2=2\norm{\Delta\bz_1}^2 + 2 (\bz_2^T\Delta\bz_1) (\bz_1^T\Delta\bz_1).
$$
Since $\bz_1$ and $\bz_2$ are unit vectors, we have:
$$
\norm{\mb{Q}_{12}-\mb{I}}
\;\leq \; 4\norm{\Delta\bz_1}.
$$
The final $(2,1)$ subblock  can be analyzed similarly. We would like to find $\mb{W}_{21}^*$ close to $\mb{W}_{21}$ such that $\mb{W}_{21}^*\mb{v}_1^*=\mb{0}$. We can define $\bar{\mb{v}}_1$, $\by_1$, and $\by_2$, and $\mb{Q}_{21}$ in a similar way to $\bar{\mb{u}}_1$, $\bz_1$, and $\bz_2$, and $\mb{Q}_{12}$. We then have $\mb{W}_{21}^*=\mb{W}_{21}\mb{Q}_{21}$ and $\mb{V}_{21}^*=(\lambda^*\mb{A}_{21}-\mb{W}_{21}^*)/\theta$.  We also obtain
$$
\Delta\by_1=\frac{1}{\sqrt{2 + \norm{\mb{v}_1^*+\bar{\mb{v}}_1}}}\Delta\bar{\mb{v}}_1,
$$
and
$$
\norm{\mb{Q}_{21}-\mb{I}}\leq 4\norm{\Delta\by_1}.
$$
Finally, we need to prove $\norm{\mb{V}^*}_{\infty}\leq 1$ and $\norm{\mb{W}}_2 \leq 1$.
By noting that $\norm{(\mb{Q}_{12}-\mb{I})\mb{W}_{12}}_{\infty}
\leq \norm{(\mb{Q}_{12}-\mb{I})\mb{W}_{12}}_2
\leq \norm{\mb{Q}_{12}-\mb{I}}_2\norm{\mb{W}_{12}}_2
\leq \norm{\mb{Q}_{12}-\mb{I}}$ and
$\norm{\Delta\bz_1} \leq 2/\sqrt{3}\eps$,
we have
$$
\ba{rl}
\norm{\mb{V}_{12}^*}_{\infty}-\norm{\mb{V}_{12}}_{\infty}&\leq
\norm{\mb{V}_{12}^*-\mb{V}_{12}}_{\infty}\\
\quad & \leq\displaystyle\theta^{-1}\left[\card{\Delta\lambda}
\norm{\mb{A}_{12}}_{\infty}+\norm{(\mb{Q}_{12}-\mb{I})
\mb{W}_{12}}_{\infty}\right]\\
 & \leq\displaystyle\theta^{-1}\left[\card{\Delta\lambda}
 \norm{\mb{A}_{12}}_{\infty}+\norm{\mb{Q}_{12}-\mb{I}}\right]\\
\quad & \leq\displaystyle\theta^{-1}\left[\card{\Delta\lambda}
\norm{\mb{A}_{12}}_{\infty}+4\norm{\Delta\bz_1}\right]
\\
 & \leq\displaystyle\theta^{-1}\left[\norm{\mb{A}_{12}}_{\infty}+5\right]\eps.
\ea
$$
Thus
$$
\norm{\mb{V}_{12}^*}_{\infty}\leq \norm{\mb{V}_{12}}_{\infty}+\theta^{-1}
\left[\norm{\mb{A}_{12}}_{\infty}+5\right]\eps\leq 1.
$$
Similarly, we also have
$$
\norm{\mb{V}_{21}^*}_{\infty}\leq \norm{\mb{V}_{21}}_{\infty}+\theta^{-1}
\left[\norm{\mb{A}_{21}}_{\infty}+5\right]\eps\leq 1.
$$
Now consider $\mb{W}^*_2$. Clearly $\norm{\mb{W}^*}_2\leq\norm{\mb{W}}_2+\norm{\mb{W}^*-\mb{W}}_2$. We have
$$
\mb{W}^*-\mb{W}=(\lambda^*-\lambda)\begin{pmatrix}
\mb{A}_{11} & \mb{0}\\
\mb{0} & \mb{A}_{22}
\end{pmatrix} + \begin{pmatrix}
\mb{0} & (\mb{Q}_{12}-\mb{I})\mb{W}_{12}\\
\mb{W}_{21}(\mb{Q}_{21}-\mb{I}) & \mb{0}
\end{pmatrix} - \begin{pmatrix}
\mb{u}_1^*(\mb{v}_1^*)^T-\mb{u}_1\mb{v}_1^T & \mb{0}\\
\mb{0} & \mb{0}
\end{pmatrix}.
$$
Thus we have:
$$
\ba{rl}
\norm{\mb{W}^*-\mb{W}}_2 &\leq \card{\Delta\lambda}\max\{\norm{\mb{A}_{11}}_2,\norm{\mb{A}_{22}}_2\}
+\max\{\norm{(\mb{Q}_{12}-\mb{I})\mb{W}_{12}}_2,
\norm{\mb{W}_{21}(\mb{Q}_{21}-\mb{I})}_2\}
\\
\quad &\quad + \norm{\mb{u}_1^*(\mb{v}_1^*)^T-\mb{u}_1\mb{v}_1^T}_2\\
\quad &\leq\max\{\norm{\mb{A}_{11}}_2,\norm{\mb{A}_{22}}_2\}\eps + 5\eps + 2\eps + \eps^2\\
\quad &\leq\displaystyle\left[\norm{\mb{A}}_2+7.5\right]\eps,
\ea
$$
since $\max\{\norm{\mb{A}_{11}}_2,\norm{\mb{A}_{22}}_2\}\leq\norm{\mb{A}}_2$ and $0<\eps<1/2$. Thus
$$
\norm{\mb{W}^*}_2\leq\norm{\mb{W}}_2+\left[\norm{\mb{A}}_2+7.5\right]\eps\leq 1.
$$
We have constructed $\mb{W}^*$ and $\mb{V}^*$ that satisfy the optimality condition for $\bX^*$, which implies that $\bX^*$ is indeed an optimal solution of Problem \refs{eq:mprob}.
\end{pf}

Proposition \ref{prop:eps} show that given an $\eps$-solution $(\lambda,\mb{u}_1,\mb{v}_1)$ of (\ref{eq:newton}), which for example, can be obtained from the current solution $(\bX_k,\bY_k)$ of the proximal point algorithm, if we could find $\mb{W}$ and $\mb{V}$ that satisfy the $\eps$-optimality conditions given in Proposition \ref{prop:eps}, then we can stop the algorithm with an accurate rank-one solution for Problem \refs{eq:mprob}.
Given $(\lambda,\mb{u}_1,\mb{v}_1)$. Let us consider the following optimization problem:
\be
\label{eq:epswv}
\ba{rl}
\min & \norm{\mb{W}}_2\\
\st & \mb{W}_{11}=(\lambda\mb{A}_{11}-\theta\mb{E}_{M\times N})-\mb{u}_1\mb{v}_1^T,\\
\quad & \mb{W}_{12}^T\mb{u}_1=\mb{0},\\
\quad & \mb{W}_{21}\mb{v}_1=\mb{0},\\
\quad & \norm{\mb{W}_{12}-\lambda\mb{A}_{12}}_{\infty}\leq\theta-
(\norm{\mb{A}_{12}}_{\infty}+5)\eps,\\
\quad & \norm{\mb{W}_{21}-\lambda\mb{A}_{21}}_{\infty}\leq\theta-
(\norm{\mb{A}_{21}}_{\infty}+5)\eps,\\
\quad & \norm{\mb{W}_{22}-\lambda\mb{A}_{22}}_{\infty}\leq\theta.
\ea
\ee

Clearly, if we could find a feasible solution of Problem \refs{eq:epswv} with the objective $\norm{\mb{W}}_2\leq 1-(\norm{\mb{A}}_2+7.5)\eps$, then the $\eps$-optimality conditions for $(\lambda,\mb{u}_1,\mb{v}_1)$ are satisfied. This problem is a non-smooth convex constrained optimization problem and our main purpose is to find a feasible solution with the objective value that is small enough. Therefore, we can simply apply the projected subgradient method
 to solve it. The projected subgradient method uses the iteration
$$
\mb{W}_{k+1}=\Pi_{{\cal W}}\left[\mb{W}_k-\alpha_k\mb{G}_k\right],
$$
where $\mb{G}_k\in\partial\norm{\mb{W}_k}_2$ is a subgradient of $\norm{.}_2$ at $\mb{W}_k$ and ${\cal W}$ is the feasible set of Problem \refs{eq:epswv}. The step size $\alpha_k$ can be chosen as one of the standard step sizes of the general subgradient method. For this problem, we choose $\displaystyle\alpha_k=O\left(1/\sqrt{k}\right)$. According to Zi\c{e}tak \cite{Zietak93}, we can always select $\mb{G}_k=\mb{u}_k\mb{v}_k^T\in\partial\norm{\mb{W}_k}_2$, where $(\mb{u}_k,\mb{v}_k)$ is the singular vectors corresponding to the largest singular value of $\mb{W}_k$. We now consider the projection problem $\Pi_{{\cal W}}(\bar{\mb{W}})$:
\be
\label{eq:proj}
\ba{rl}
\Pi_{{\cal W}}(\bar{\mb{W}})\in\arg\min & \norm{\mb{W}-\bar{\mb{W}}}_F^2\\
\st & \mb{W}_{11}=(\lambda\mb{A}_{11}-\theta\mb{E}_{M\times N})-\mb{u}_1\mb{v}_1^T,\\
\quad & \mb{W}_{12}^T\mb{u}_1=\mb{0},\\
\quad & \mb{W}_{21}\mb{v}_1=\mb{0},\\
\quad & \norm{\mb{W}_{12}-\lambda\mb{A}_{12}}_{\infty}\leq\theta-
(\norm{\mb{A}_{12}}_{\infty}+5)\eps,\\
\quad & \norm{\mb{W}_{21}-\lambda\mb{A}_{21}}_{\infty}\leq\theta
-(\norm{\mb{A}_{21}}_{\infty}+5)\eps,\\
\quad & \norm{\mb{W}_{22}-\lambda\mb{A}_{22}}_{\infty}\leq\theta.
\ea
\ee

The objective function $\norm{\mb{W}-\bar{\mb{W}}}_F^2$ is element-wise separable; therefore, Problem \refs{eq:proj} is block-wise separable. For the
$(1,1)$ subblock, we have the fixed solution $\mb{W}_{11}=(\lambda\mb{A}_{11}-\theta\mb{E}_{M\times N})-\mb{u}_1
\mb{v}_1^T$. For the $(2,2)$ subblock, it is a simple element-wise separable optimization problem:
$$
\ba{rl}
\min & \norm{\mb{W}_{22}-\bar{\mb{W}}_{22}}_F^2\\
\st & \norm{\mb{W}_{22}-\lambda\mb{A}_{22}}_{\infty}\leq\theta,
\ea
$$
whose optimal solution can be computed as follows:
$$
\mb{W}_{22}=\max\left\{\min\left\{\bar{\mb{W}}_{22},\lambda\mb{A}_{22}+\theta\right\},\lambda\mb{A}_{22}-\theta\right\}.
$$
For the $(1,2)$ subblock, the corresponding optimization problem is column-wise separable:
$$
\ba{rl}
\min & \norm{\mb{W}_{12}-\bar{\mb{W}}_{12}}_F^2\\
\st & \mb{W}_{12}^T\mb{u}_1=\mb{0},\\
\quad & \norm{\mb{W}_{12}-\lambda\mb{A}_{12}}_{\infty}
\leq\theta-(\norm{\mb{A}_{12}}_{\infty}+5)\eps.
\ea
$$
Each subproblem is a quadratic knapsack problem which can be written as follows:
$$
\ba{rl}
\min & \norm{\mb{w}-\bar{\mb{w}}}^2\\
\st & \mb{u}_1^T\mb{w}=0,\\
\st & \mb{l}\leq\mb{w}\leq\mb{u}.
\ea
$$
According to Brucker \cite{Brucker84}, there is an $O(n)$ algorithm for these quadratic knapsack problems. Thus we can find $\mb{W}_{12}$ efficiently. Similarly, the $(2,1)$ subblock  can be found by solving a number of quadratic knapsack problems since the corresponding optimization problem for it is row-wise separable.

\section{Numerical Examples}
\subsection{Sailboat Bitmap Image Example}
\begin{figure}[htp]
\centering
\includegraphics[width=0.7\textwidth]{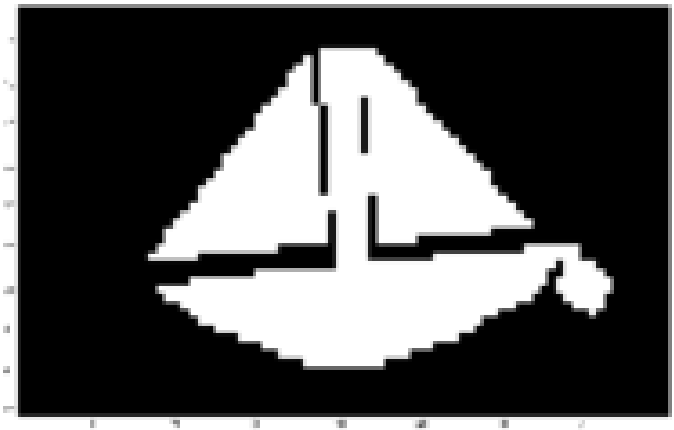}
\caption{Bitmap image of the complete sailboat}
\label{fig:sailboat}
\end{figure}

\begin{figure}[htp]
\centering
\includegraphics[width=0.9\textwidth]{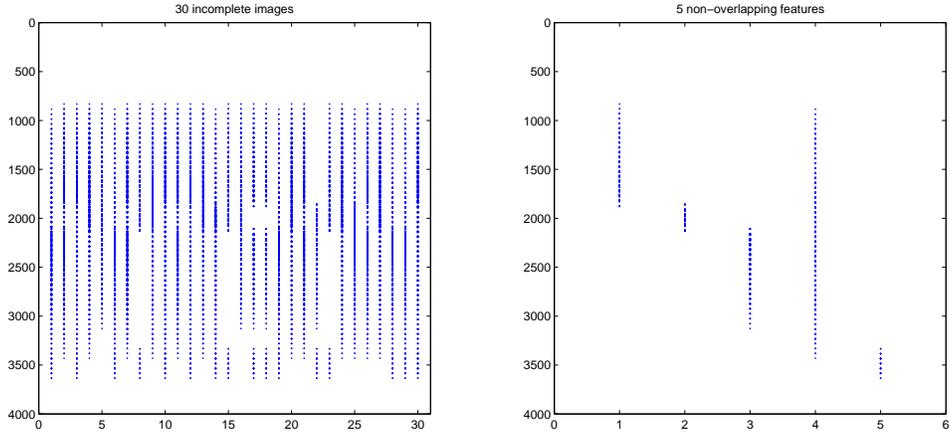}
\caption{Collection of $30$ incomplete images of the sailboat and its $5$ components}
\label{fig:imgcol}
\end{figure}

In this example, we use a $80$-by-$50$ black-and-white bitmap image of a sailboat. There are $5$ distinct components or non-overlapping features in this image: left sail (Feature 1), sail mast (Feature 2), right sail (Feature 3), hull (Feature 4), and rudder (Feature 5). The bitmap image of the sailboat is shown in Figure \ref{fig:sailboat}. A matrix $\mb{A}$ is created with $30$ columns, each of which represents a bitmap image of the sailboat with just $3$ out of $5$ features. The matrix $\mb{A}$ therefore has $5$ rank-one submatrices composed of all ones since the bitmap image is black-and-white. The structure of matrix $\mb{A}$ and all the features (in terms of non-zero elements) are shown in Figure \ref{fig:imgcol}. We would like to use our proposed formulation to extract these rank-one submatrices.

Our main task in this example is to use our proposed formulation to extract the features from the matrix $\mb{A}$. We have developed two algorithms to solve Problem \refs{eq:mprob}, the primal and dual. For these numerical examples, we will use the dual algorithm mainly due to its accuracy with respect to optimal solutions. This superior accuracy could be explained by the fact that the subproblem we solved in each iteration of the dual algorithm is similar to the original problem. The stopping criterion stated in Proposition \ref{prop:eps} is implemented. We test the conditions of Kantorovich's theorem, solve the system of equations in \refs{eq:newton} to obtain an $\eps$-solution, and then use the projected subgradient method to find a feasible solution $\mb{W}$ of Problem \refs{eq:epswv}. Since the projection problem with $\bar{\mb{W}}=\mb{0}$ can be considered as a relaxation of Problem \refs{eq:epswv} in which the spectral norm is replaced by the Frobenius norm, we will start the projected subgradient method with $\mb{W}_0=\mb{0}$. The stopping criteria for this method are the condition $\norm{\mb{W}}_2\leq 1-(\norm{\mb{A}}_2+7.5)\eps$, the maximum number of iterations, and the change in objective values. As the testing process is computationally quite expensive, therefore we only use it once per a fixed number (say, 10) of outer iterations.

For each value of the parameter $\theta$, we will obtain the optimal solution $\bX_1^*$ and $\bX_2^*$. The decision variable $\bX_2$ corresponds to the $l_1$-norm part of the objective function; therefore, we use the sparsity structure of $\bX_2^*$ to construct the final solution $\bX_F^*$ with the elements of $\bX_2^*$. According to Theorem 5 in \cite{Doan10}, the optimal solution of Problem \refs{eq:mprob} indicates the exact sparsity structure of the rank-one submatrix (even under small random noise) with appropriate $\theta$. Therefore, in this experiment, we extract the rank-one approximation of $\bX_F^*$ and use its sparsity structure as the sparsity structure of the extracted feature. Next, we need to select an appropriate value for $\theta$. For $\theta\approx 0$, the algorithm returns the rank-one approximation of matrix $\mb{A}$, which is an average of all features and for the purpose of extracting single features, this averaging effect is not desirable. On the other hand, we prefer large submatrices (large features) over small ones. Similar to the L-curve method used to select a regularization parameter, we construct the curve ${\cal L}:=\{\norm{\bX(M_\theta,N_\theta)}_F,\norm{\mb{A}(M_\theta,N_\theta)-\bX(M_\theta,N_\theta)}_F,\,\theta\geq 0\}$, where $(M_{\theta},N_{\theta})$ is the sparsity structure obtained from the algorithm and $\bX(M_{\theta},N_{\theta})$ is the rank-one approximation of $\mb{A}(M_\theta,N_\theta)$. We then pick $\theta$ that balances the feature largeness, $\norm{\bX(M_\theta,N_\theta)}_F$, and feature averaging measure $\norm{\mb{A}(M_\theta,N_\theta)-\bX(M_\theta,N_\theta)}_F$. After selecting $\theta$, we obtain $\bX(M_{\theta},N_{\theta})=\mb{u}(M_{\theta})\mb{v}(N_{\theta})^T$, where $\max(\mb{v}(N_{\theta}))=1$. The vector $\mb{u}(M_{\theta})$ represents the extracted feature and $\mb{v}(N_{\theta})$ indicates how significant the feature is in each boat image. After extracting a feature, we remove the feature from the image by setting $\mb{A}(M_{\theta},N_{\theta})=\mb{0}$ and continue to find new (non-overlapping) features. This method for choosing $\theta$ is clearly just a heuristic and a more concrete approach for $\theta$ selection is still an important issue for future research.

We are now ready to run our algorithm on this sailboat example. We set the main tolerance to be $\eps=10^{-6}$, the maximum number of iterations to be $1000$, and for each subproblem, the maximum number of iterations is set to be $30$. There is also the parameter $\lambda$ of the proximal point framework that we need to select. This parameter controls the convergence of the algorithm and for this example, $\lambda=O(1/\theta)$ works well most of the time. We can always adjust $\lambda$ (and number of iterations) to get better convergence if the initial setting does not achieve the tolerance required. We set $\eps_s=10^{-10}$ as the tolerance used in Newton's method to test the additional stopping criterion. Finally, the values of $\theta$ are selected uniformly from three different ranges, small range $[0.01,0.1]$, medium range $[0.1,1]$, and large range $[1,10]$, $10$ values in each range.

We start with the matrix $\mb{A}$. Except for the first value of $\theta$ ($\theta=0.01$), all other values result in the same rank-one submatrix of $\mb{A}$, which means $\norm{\mb{A}(M_{\theta},N_{\theta})-\bX(M_{\theta},N_{\theta})}=0$. Thus we do not need to use the curve $\cal L$ and just need to pick any value of $\theta>0.01$. The vector $\mb{v}({N_\theta})$ is a zero-one vector indicating that either the feature $\mb{u}(M_{\theta})$ appears completely in an image or it does not appear at all. The feature $\mb{u}(M_{\theta})$ represents the exact combination of Feature 1 and 4, which is the left sail and the hull.

\begin{figure}[htp]
\centering
\includegraphics[width=0.7\textwidth]{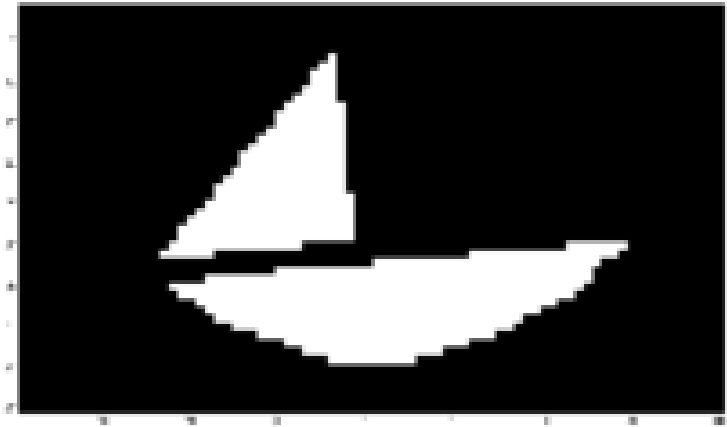}
\caption{First extracted feature: the combination of left sail and hull}
\label{fig:firstfeat}
\end{figure}

We now exclude the first extracted feature from all the images and continue to find new (non-overlapping) features. Table \ref{table:features} shows all the features that we obtain with the size of the features ($s_i$), number of images that share each feature ($n_i$), and their description.

\begin{table}[htp]
\begin{center}
\begin{tabular}{|c|c|c|l|}
\hline
$i$ & $s_i$ & $n_i$ & Description \\
\hline
$1$ & $710$ & $15$ & Left sail and hull\\
$2$ & $694$ & $6$ & Right sail and hull\\
$3$ & $252$ & $10$ & Right sail\\
$4$ & $156$ & $5$ & Sail mast and rudder\\
$5$ & $268$ & $7$ & Left sail\\
$6$ & $119$ & $9$ & Sail mast\\
$7$ & $439$ & $1$ & Hull\\
$8$ & $34$ & $11$ & Rudder\\
\hline
\end{tabular}
\caption{All extracted features obtained from the algorithm}
\label{table:features}
\end{center}
\end{table}

The results show that our algorithm can pick out the large common features that are inherent in the structure of the image set. For example, a combination of our defined features is indeed a large common feature if there are enough images that share that combination of features.

To end this section, we would like to comment on the efficiency of the additional stopping criterion based on Proposition \ref{prop:eps}. When the test indicates the convergence is achieved, it is guaranteed that the supports of $\mb{u}$ and $\mb{v}$ have been correctly identified for the rank-one optimal solution $\bX$. On the other hand, because the test uses heuristics to find multipliers, it may occur that the optimal supports are attained and yet the test fails to indicate that. In this sailboat example, we ran the algorithm with $8$ different matrices, $\mb{A}_0=\mb{A}$, $\mb{A}_1,\ldots,\mb{A}_7$ with subsequent extraction of features one by one. The additional stopping criterion works for $4$ out of $8$ matrices and we obtain a significant reduction in both computational time and number of iterations while maintaining highly accurate solutions ($\eps_s=10^{-10}$). Table \ref{table:improv} shows these improvements with $\theta=0.2$, where (DDPA)/(ADDPA) is the algorithm without/with the additional stopping criteria. 
\begin{comment} The tuple includes the number of outer iterations, the total number of inner iterations, the total computational time, and the time allocated for testing the additional stopping criterion.\end{comment} 
The number of iterations with (ADDPA) is either $10$ or $20$ since in this example, we only test the additional stopping criterion once per $10$ (outer) iterations. We can see that there are cases when the additional stopping criterion can be used very early to stop the algorithm with a guaranteed highly accurate solution.

\begin{table}[htp]
\begin{center}
\begin{tabular}{|c|c|c|}
\hline
Matrix & (DPPA)  &  (ADPPA)\\
\hline
$\mb{A}_4$ & $(59,325,9.36s,0.00s)$ & $(20,123,6.36s,2.32s)$\\
$\mb{A}_5$ & $(52,267,7.28s,0.00s)$ & $(10, 60,4.23s,2.39s)$\\
$\mb{A}_6$ & $(62,499,12.4s,0.00s)$ & $(20,187,6.96s,1.84s)$\\
$\mb{A}_7$ & $(29,148,3.59s,0.00s)$ & $(10, 57,4.00s,2.48s)$\\
\hline
\end{tabular}
\caption{Outer iteration number, inner iteration number, total computational time, and convergence testing time for (DDPA)/(ADDPA)}
\label{table:improv}
\end{center}
\end{table}

\subsection{Image Database Test Case}
We conduct the experiment on the Frey face dataset, which consists of $1965$ registered face images of size $28\times 20$. The matrix $\mb{A}$ has the size of $1965\times 560$, where each column represents a single face image. We again use the dual algorithm with the additional stopping criterion and maintain all parameters the same as in the previous example. The additional stopping criterion is less effective in this test case. However, when it works, we again have a significant improvement in computational time and number of iterations. For example, with $\mb{A}_0=\mb{A}$ and $\theta=0.2$, we have the following results for (DDPA) and (ADDPA) respectively: $(245,6466,1.21\times 10^3s,0.00s)$ and $(100,2140,4.25\times 10^2s,23.7s)$, where the tuple is explained in the caption of Figure \ref{table:improv}.

\begin{figure}[htp]
\centering
\includegraphics[width=1.0\textwidth]{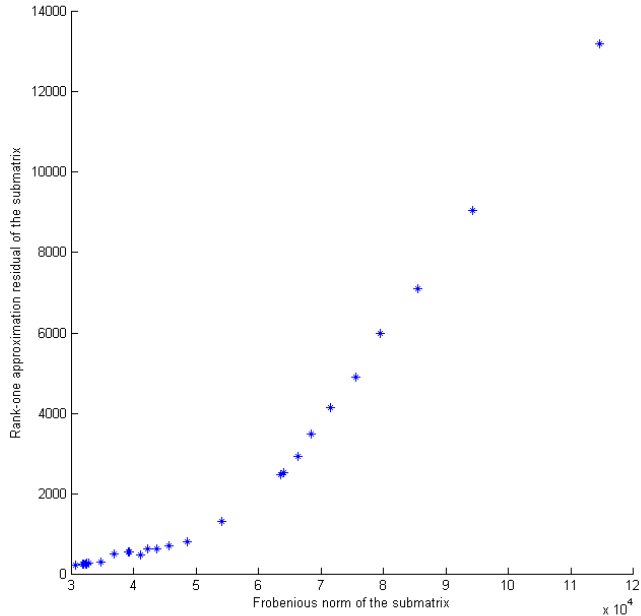}
\caption{Feature largeness vs. feature averaging measure for different $\theta$}
\label{fig:measures}
\end{figure}

\begin{figure}[htp]
\centering
\includegraphics[width=1.0\textwidth]{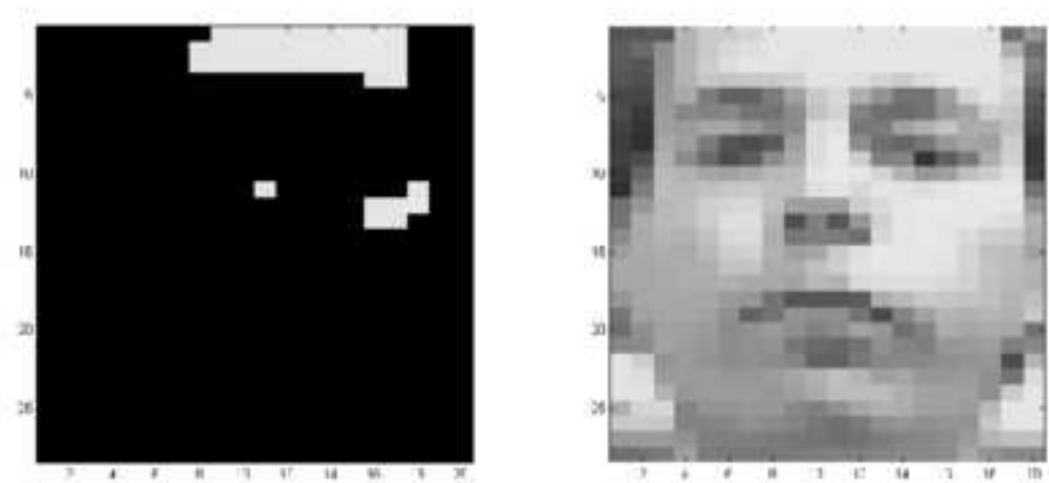}
\caption{First feature and the image that has the highest significance factor}
\label{fig:facefeat}
\end{figure}

\begin{figure}[htp]
\centering
\includegraphics[width=1.0\textwidth]{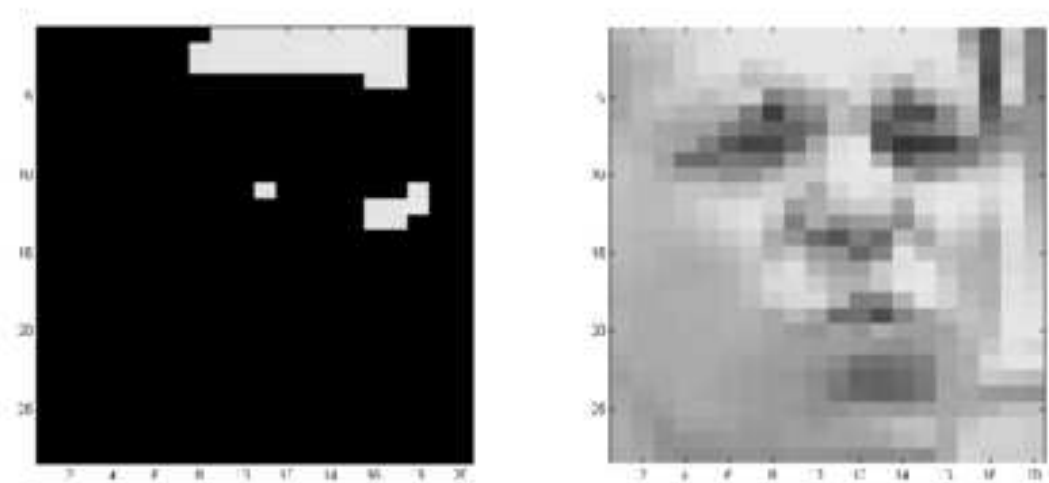}
\caption{An image without the first feature}
\label{fig:nofeat}
\end{figure}

We apply the algorithm to the matrix $\mb{A}$ and Figure \ref{fig:measures} shows the curve $\cal L$ obtained with different values of $\theta$. We select $\theta=0.2$ at the curviest point on $\cal L$, which indicates the balance between feature largeness and feature averaging measure. We obtain the feature $\mb{u}_1$ and the significance vector $\mb{v}_1$ indicating how strong the appearance of that feature is in each image. The feature is composed of $38$ pixels and there are $1557$ images that are considered to have this feature with the significance factor of at least $95.14\%$, where the significance factor of the feature in image $j$ is defined as $v_1(j)/\norm{\mb{v}_1}_{\infty}$. Figure \ref{fig:facefeat} presents the first feature and the face image that has the significance factor of $100\%$ for this feature. Basically, the first feature shows the right forehead, a part of right cheekbone, and the tip of the nose. This feature is common among the images ($1557$ of them), there are images that do not share the feature. Figure \ref{fig:nofeat} shows one of such images.

We remove the first feature from all images that share that feature and continue to find new features. Table \ref{table:feats} shows the size of the feature $i$ ($s_i$), number of images that share the feature $i$ ($n_i$), and the minimum significance factor for each feature $i$ ($f_i^{\min}$), $i=1,\ldots,10$.

\begin{table}[htp]
\begin{center}
\begin{tabular}{|c|c|c|c||c|c|c|c|}
\hline
$i$ & $s_i$ & $n_i$ & $f_i^{\min}$ &$i$ & $s_i$ & $n_i$ & $f_i^{\min}$\\
\hline
$1$ & $38$ & $1557$ & $95.14\%$ & $6$ & $28$ & $673$ & $83.27\%$\\
$2$ & $27$ & $896$ & $92.19\%$ & $7$ & $21$ & $578$ & $80.38\%$\\
$3$ & $29$ & $1096$ & $87.61\%$ & $8$ & $20$ & $555$ & $87.59\%$\\
$4$ & $24$ & $847$ & $83.12\%$ & $9$ & $35$ & $291$ & $80.73\%$\\
$5$ & $25$ & $791$ & $83.67\%$ & $10$ & $13$ & $598$ & $71.31\%$\\
\hline
\end{tabular}
\caption{Information of the first ten extracted features}
\label{table:feats}
\end{center}
\end{table}

\begin{figure}[htp]
\centering
\includegraphics[width=1.0\textwidth]{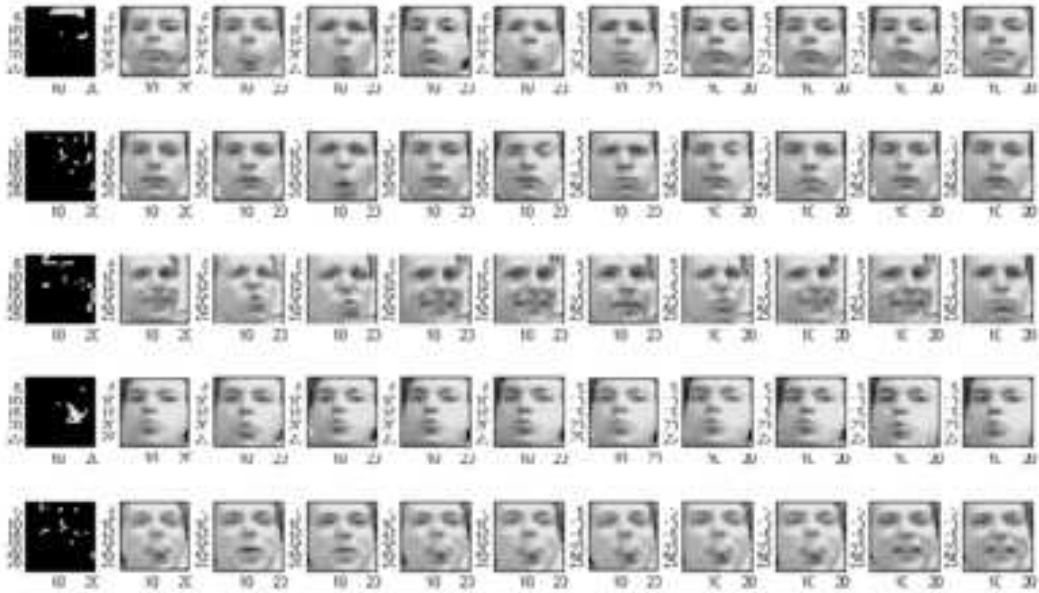}
\caption{First five features and images with highest significance factors}
\label{fig:feats15}
\end{figure}

\begin{figure}[htp]
\centering
\includegraphics[width=1.0\textwidth]{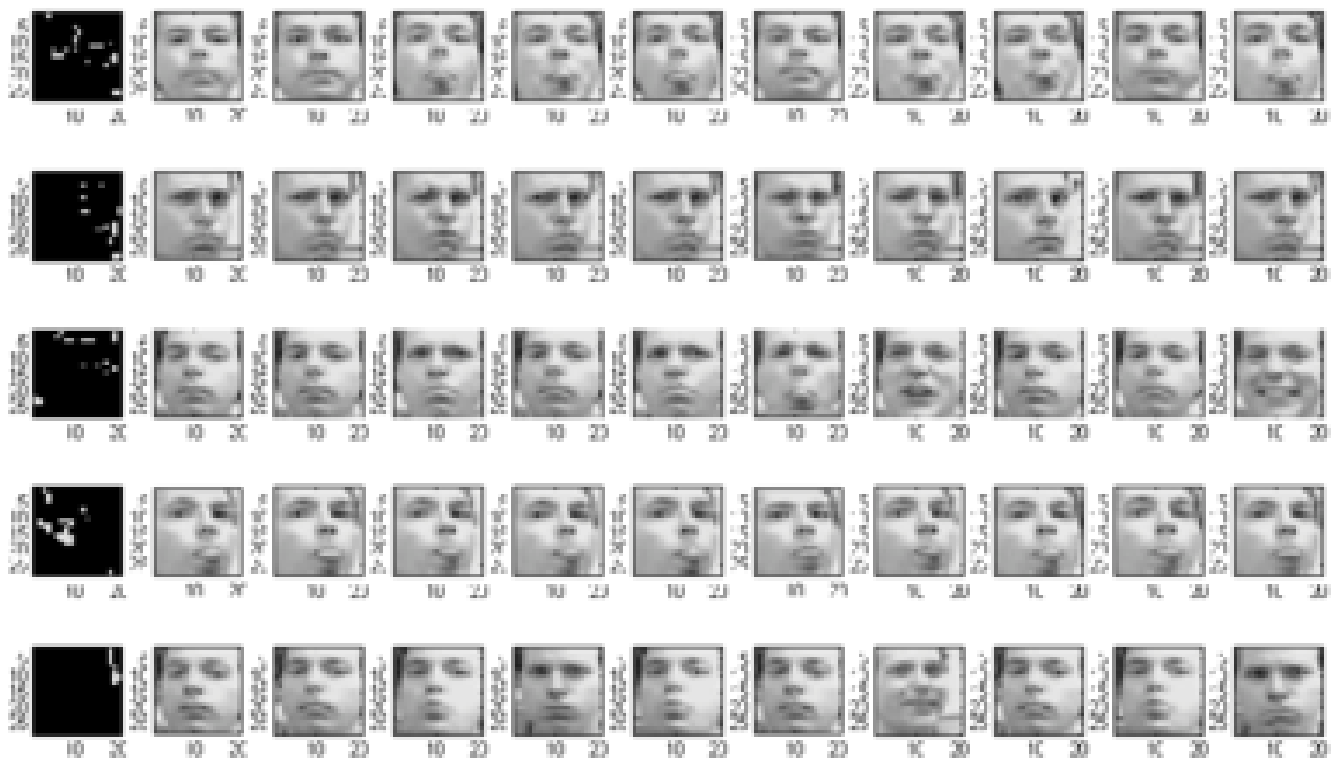}
\caption{Feature 6 to Feature 10 and images with highest significant factors}
\label{fig:feats610}
\end{figure}

Figure \ref{fig:feats15} and \ref{fig:feats610} presents each feature and the ten face images that have the highest significance factors for that feature.

\begin{figure}[htp]
\centering
\includegraphics[width=1.0\textwidth]{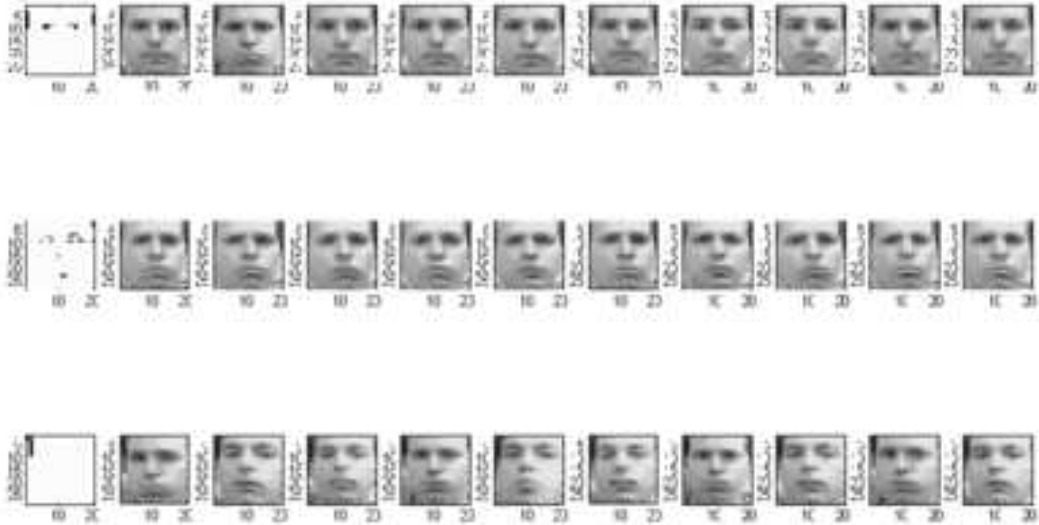}
\caption{First three negative features and images with highest significance factors}
\label{fig:nfeats13}
\end{figure}

The features are not easy to observe or distinguish. However, with images that have high significance factors, we can see that some features could be associated with a certain orientation of the face or lighting of images. For example, Feature 4 and Feature 9 clearly show the right (or left) cheek when Frey faces left (or right). Certain lighting of the background can also define features, which is the case of Feature 3 and Feature 7. Another observation is that since this approach favors large submatrices, which means other features defined by small entries (dark pixels), for examples, eyes or mouth, will not be picked up as major features. In this particular application of visual features, we can define \emph{negative} features, which correspond to the features of the negative images. In order to find these negative features, we construct the negative images and apply the algorithm to this set of images. In this example, the algorithm is applied to $\mb{B}=255\mb{E}-\mb{A}$, where $\mb{E}$ is the matrix of all ones. The coefficient $255$ appears due to the range of pixel intensities in these images. For each feature $\mb{u}$ extracted from $\mb{B}$, we define $\mb{u}_n=255\mb{e}-\mb{u}$, where $\mb{e}$ is the vector of all ones, as the negative feature of the original set of images. The first three extracted features are presented in Figure \ref{fig:nfeats13}. The first negative feature has both straight dark eyes with two dark background columns at both sides on top. The second one focuses on the darker right eye, the left nostril, and also the chin. And the third one is a long dark background column on the top left.

\newpage

\renewcommand{\baselinestretch}{1.00}
\small
\bibliographystyle{plain}
\bibliography{NSRFFF}

\begin{thebibliography}{1}

\bibitem{BT09}
A.~Beck and M.~Teboulle.
\newblock A fast iterative shrinkage-thresholding algorithm for linear inverse
  problems.
\newblock {\em SIAM Journal on Imaging Sciences}, 1:183--202, 2009.

\bibitem{Brucker84}
P.~Brucker.
\newblock An ${O}(n)$ algorithm for quadratic knapsack problems.
\newblock {\em Operations Research Letters}, 3:163--166, 1984.

\bibitem{Doan10}
X.~V. Doan and S.~Vavasis.
\newblock Finding approximately rank-one submatrix.
\newblock Under review, SIAM Journal of Optimization, URL:
  \url{http://arxiv.org/abs/1011.1839}, 2010.

\bibitem{Liu09}
Y.~J. Liu, D.~Sun, and K.~C. Toh.
\newblock An implementable proximal point algorithmic framework for nuclear
  norm minimization, 2009.
\newblock Preprint.

\bibitem{Rockafellar70}
R.~T. Rockafellar.
\newblock {\em Convex Analysis}.
\newblock Princeton University Press, Princeton, NJ, 1970.

\bibitem{Rockafellar76}
R.~T. Rockafellar.
\newblock Monotone operators and the proximal point algorithm.
\newblock {\em SIAM Journal of Control and Optimization}, 14:877--898, 1976.

\bibitem{Tapia71}
R.~A. Tapia.
\newblock The {K}antorovich theorem for {N}ewton method.
\newblock {\em The American Mathematical Monthly}, 78:389--392, 1971.

\bibitem{Toh09}
K.~C. Toh and S.~W. Yun.
\newblock An accelerated proximal gradient algorithm for nuclear norm
  regularized least squares problems, 2009.
\newblock To appear in Pacific Journal of Optimization.

\bibitem{Zietak93}
K.~Zi\c{e}tak.
\newblock Properties of linear approximations of matrices in the spectral norm.
\newblock {\em Linear Algebra Applications}, 183:41--60, 1993.

\end{thebibliography}
\end{document}